\documentclass[12pt,centertags,oneside]{amsart}
\usepackage{amsmath,amstext,amsthm,amscd,typearea}
\usepackage{amssymb}
\usepackage{a4wide}
\usepackage[mathscr]{eucal}
\usepackage{mathrsfs}
\usepackage{typearea}
\usepackage{charter}
\usepackage{pdfsync}
\usepackage{url}
\usepackage{xcolor}

\usepackage[a4paper,width=16.2cm,top=3cm,bottom=3cm]{geometry}

\numberwithin{equation}{section}

\newtheorem{theorem}{Theorem}[section]
\newtheorem{definition}[theorem]{Definition}
\newtheorem{proposition}[theorem]{Proposition}

\newtheorem{lemma}[theorem]{Lemma}
\newtheorem{remark}[theorem]{Remark}

\newcommand{\cali}[1]{\mathscr{#1}}

\newcommand{\Aut}{{\rm Aut}}

\newcommand{\dist}{\mathop{\mathrm{dist}}\nolimits}

\renewcommand{\Im}{\mathop{\mathrm{Im}}\nolimits}

\newcommand{\ddc}{\text{\normalfont dd}^c}

\newcommand{\Reg}{\text{\normalfont Reg}}

\newcommand{\PSH}{{\rm PSH}}

\newcommand{\id}{{\rm id}}

\newcommand{\codim}{{\rm codim\ \!}}

\newcommand{\FS}{{\rm FS}}
\newcommand{\B}{\mathbb{B}}
\newcommand{\C}{\mathbb{C}}

\newcommand{\N}{\mathbb{N}}
\newcommand{\Z}{\mathbb{Z}}
\newcommand{\R}{\mathbb{R}}
\renewcommand\P{\mathbb{P}}

\renewcommand{\S}{\mathbb{S}}

\title[]{Equilibrium measures of meromorphic self-maps on non-K\"ahler manifolds}

\author{Duc-Viet Vu}
\address{University of Cologne, Mathematical Institute, Germany}
\email{vuduc@math.uni-koeln.de}

\date{\today}
\begin{document}

\begin{abstract} Let $X$ be a compact complex non-K\"ahler manifold and $f$ a dominant meromorphic self-map of $X$. Examples of such maps are self-maps of Hopf manifolds, Calabi-Eckmann manifolds, non-tori nilmanifolds, and their blowups.  We prove that if $f$ has a dominant topological degree, then $f$ possesses an equilibrium measure $\mu$  satisfying well-known properties as in the K\"ahler case.   The key ingredients are the notion of  weakly d.s.h. functions substituting d.s.h. functions in the K\"ahler case and the use of suitable test functions in Sobolev spaces.  A large enough class of holomorphic self-maps with a dominant topological degree on Hopf manifolds is also given. 
\end{abstract}

\maketitle

\medskip

\noindent
{\bf Classification AMS 2010}:  32U40, 32H50, 37F05. 

\medskip



\section{Introduction} \label{introduction}

 Let $X$ be a compact complex  manifold of dimension $k$. Let $f$ be a  dominant meromorphic self-map  of $X.$   Let $\omega$ be a  strictly positive Hermitian $(1,1)$-form on $X.$  For $0 \le q \le k,$ put
$$d_q(f):= \limsup_{n \to \infty}\bigg( \int_X  (f^n)^* \omega^q \wedge  \omega^{k-q}\bigg)^{1/n}.$$
We will write $d_q$ for $d_q(f)$ if no confusion arises.   We can see easily that $d_q$ is independent of the choice of $\omega.$ The number $d_0$ is always $1$ and $d_k$  is  the topological degree of $f.$ When $f$ is holomorphic, $d_q$ is finite because  the differential of $f$ is of $L^\infty$-norm uniformly bounded on $X.$ We call $d_q$ the $q^{th}$ \emph{dynamical degree} of $f$ for $0 \le q \le k.$ 

When $X$ is K\"ahler, the numbers $d_q$ are crucial \emph{finite} bi-meromorphic invariants of $f;$ see  \cite{DS_upperbound,Ds_upperbound_mero,DS_regula}.   We don't know whether $d_q$ for $1\le q \le k-1$  is finite for general  $X.$  In what follows, we will study the dynamics of $f$ with  $d_k(f)> d_{k-1}(f).$  In the K\"ahler case, such a map is said to have a dominant topological degree and its dynamics has been thoroughly investigated, see \cite{DS_book} and references therein for information.  We emphasize that in our context,  it is not clear whether the assumption $d_k> d_{k-1}$ implies $d_k> d_q$ for $1 \le q\le k-1$ as in the K\"ahler case.  

A \emph{quasi-p.s.h.} function on $X$ is a function from $X$ to $[-\infty, \infty)$ which is locally the sum of a plurisubharmonic function and a smooth one.  For a given continuous $(1,1)$-form $\eta,$ denote by $\PSH_0(\eta)$ the set of quasi-p.s.h. functions $\varphi$ such that $\ddc \varphi +\eta \ge 0$ and $\sup_X \varphi=0.$ Equip $\PSH_0(\eta)$ with the induced distance from $L^1(X)$ by using  the natural  inclusion $\PSH_0(\eta) \subset L^1(X).$ 

Recall from \cite{DS_tm} that  a complex measure $\mu$ on $X$ is said to be \emph{PC} if every quasi-p.s.h. function is $\mu$-integrable  and for every sequence $(\varphi_n)_{n \in \N}$ of quasi-p.s.h. functions converging to $\varphi$ in $L^1$  such that  $\ddc \varphi_n+ \eta \ge 0$ for some  smooth form $\eta$ independent of $n,$  we have $\langle \mu, \varphi_n \rangle \to \langle \mu, \varphi \rangle$.

 A \emph{pluripolar set} in $X$ is a subset of $X$ contained in $\{\varphi=-\infty\}$ for some quasi-p.s.h. function $\varphi.$ By \cite{Vu_pluripolar},  every locally pluripolar set in $X$ is pluripolar. This result implies in particular  that there exist abundantly singular quasi-p.s.h. functions on $X.$ Observe that every PC measure  has no mass on pluripolar sets.     Here is our  first main result.
 
\begin{theorem}\label{the_topodereedi} Let $X$ be a compact complex manifold of dimension $k$ and $f$ a dominant meromorphic self-map of $X$ with $d_k> d_{k-1}$.  Let $\nu$ be a complex measure  with $L^{2k+1}$ density on $X$ so that  $\nu(X)=1$. Then  $d_k^{-n} (f^n)^*\nu$ converges weakly to a PC probability measure $\mu_f$ of entropy $\ge \log d_k$ independent of $\nu$ as $n \to \infty$ such that $d_k^{-1} f^* \mu_f= \mu_f$ and   if $f$ is holomorphic then  for every Hermitian metric $\omega$ on $X,$  $\mu_f$ is H\"older continuous on $\PSH_0(\omega).$    
\end{theorem}

The  above measure $\mu_f$ is called \emph{the equilibrium measure} of $f.$ We emphasize that unlike the K\"ahler case, it is not clear to us whether the entropy of $\mu_f$ is equal to $\log d_k.$  As to the lower bound on the entropy of $\mu_f,$ without the assumption that  $d_k>d_{k-1},$  De Th\'elin and Vigny considered a more general sequence of measures constructed from $f^n$ and proved a lower bound for the entropy of limit measures of that sequence, see Theorem 1 and the remark following it in \cite{DeThelin-Vigny} for details. 

The H\"older continuity of $\mu_f$ on $\PSH_0(\omega)$ for $f$ holomorphic implies that $\mu_f$ is \emph{moderate} in the sense that  there exist constants $\epsilon,M>0$ such that  for every $\varphi \in \PSH_0(\omega),$ we have 
$$\int_X e^{- \epsilon \varphi} d\mu_f \le M,$$
see \cite{DinhVietanhMongeampere} for a proof.  A large  class of holomorphic endomorphisms of Hopf manifolds having dominant topological degree is given in  Lemma \ref{le_ex_dominanthopf} in  Section \ref{example} and the comment following it.  
 
The existence of $\mu_f$ is proved by Fornaess-Sibony and  Russakovskii-Shiffman \cite{Sibony_pano,FS_II,Fornaess_Sibony_1994,RS} for $X=\P^k$,  Guedj \cite{Guedj} for $X$ projective and he also shows that quasi-p.s.h. functions are $\mu_f$-integrable, see also \cite{DinhSibony_allure} for the case of polynomial-like maps. The stronger fact that for $X$ K\"ahler,  $\mu_f$  is PC is    proved by   Dinh-Sibony \cite{DS_tm} by using a key property that  the space of d.s.h. functions (differences of two quasi-p.s.h. functions) is preserved by meromorphic maps.  However,  it seems that this property no longer  holds in the non-K\"ahler case. We refer  to \cite{DS_book,Duval_Briend_entropy,Lyubich,FLM} and references therein for more information on the K\"ahler case. 

In order to prove that $\mu_f$ is PC in Theorem \ref{the_topodereedi},  we  introduce a new class of functions called \emph{weakly d.s.h.} functions which, to some extent,  replace the role of d.s.h. functions (differences of two quasi-p.s.h. functions) in K\"ahler case.  These functions enjoy a compactness property similar to that of d.s.h. functions and the pull-back of d.s.h. functions by meromorphic maps are weakly d.s.h.. We also obtain the exponentially mixing property of $\mu_f$ generalizing similar results in the K\"ahler case by Dinh-Sibony in  \cite[The. 1.1]{DS_decay} and \cite[The. 1.3]{DS_tm}.  

\begin{theorem} \label{the_mixing} Let $X,f,d_{k-1}, d_k,\mu_f$ be as in Theorem \ref{the_topodereedi}. Then $\mu_f$ is exponentially mixing  in the sense that for every constant $\epsilon>0$ with $d_k> d_{k-1}+\epsilon$ and $0<\alpha \le 1,$ there exists a constant $c_{\epsilon,\alpha}$ such that   
$$\big |  \langle \mu_f, (\psi \circ f^n) \varphi\rangle - \langle \mu_f, \psi \rangle \langle \mu, \varphi \rangle  \big | \le c_{\epsilon, \alpha} \|\psi\|_\infty \|\varphi\|_{C^\alpha} (d_{k-1}+\epsilon)^{n \alpha/2} d_k^{-n \alpha/2}$$
for every $n\ge 0,$ every $\psi \in L^\infty(X)$ and every H\"older continuous function $\varphi$ of order $\alpha.$ In particular, $\mu_f$ is K-mixing.

If a real-valued H\"older continuous function $\varphi$ is not a coboundary, \emph{i.e}, there doesn't exist $\psi \in L^2(X)$ with $\varphi= \psi \circ f - \psi,$ and satisfies $\langle \mu, \varphi \rangle=0,$ then $\mu_f$ satisfies the central limit theorem, that means there is a constant $\sigma>0$ such that for every interval $I \subset \R,$ we have 
$$\lim_{n \to \infty} \mu_f \bigg \{ \frac{1}{\sqrt{n}} \sum_{j=0}^{n-1} \varphi \circ f^j \in I     \bigg\} = \frac{1}{\sqrt{2\pi} \sigma} \int_I e^{-x^2/(2 \sigma^2)} dx.$$
\end{theorem}

In the above statement, only the decay of correlation requires new arguments. The K-mixing property and the central limit theorem are deduced by using similar arguments from \cite{DS_decay,DS_book}.  Due to the same reason with the pull-back of d.s.h. functions presented above,  the arguments  in  \cite{DS_tm} couldn't be applied directly to obtain the expected decay of correlation. Our approach is based on ideas from \cite{DS_decay}: use a suitable class of functions in the Sobolev space $W^{1,2}$ as test functions.  

The above results for meromorphic maps still hold for meromorphic correspondences. But in order to keep the presentation as simple as possible, we don't elaborate it here.  In the next section, we prove Theorem \ref{the_topodereedi}. A proof of Theorem \ref{the_mixing} is given in Section \ref{secW12}.   Examples of dynamical systems on non-K\"ahler manifolds are given in Section \ref{example}. 
\\

\noindent
\textbf{Acknowledgments.} The author would like to express his gratitude to Tien-Cuong Dinh, Vi\^et-Anh Nguy\^en and Gabriel Vigny  for fruitful discussions.  This research is supported by a postdoctoral fellowship of Alexander von Humboldt Foundation.

\section{Maps with dominant topological degrees} \label{sec_mapdominant}

In this section, we will prove Theorem \ref{the_topodereedi}. For a current $T$ of order $0$ defined on a manifold $X,$ we denote by $\|T\|_X$ the mass of $T$ over $X.$  We will write $\lesssim$ (resp. $\gtrsim$) for $\le$ (resp. $\ge$) modulo a multiplicative constant independent of involving terms in the inequality.

Let $\B_r$ be the ball centered at $0$  of  radius $r$ of $\C^k,$ where $r\in \R^+.$ For $r:=1$ we put $\B:=\B_1.$  The following is crucial for us.

\begin{lemma} \label{le-11current}
Let $r\in (0,1).$ Then  for every  real closed $(1,1)$-current  $R$ of order $0$ defined on $\B,$ there is a function $U_R$ in $L^{1+1/(2k)}(\B_r)$ such that the following three properties hold:

$(i)$ $$R=\ddc U_R$$
 on $\B_r$,

$(ii)$  
$$\|U_R\|_{L^{1+1/(2k)}(\B_r)}\le c_r \|R\|_{\B}$$
for some constant $c_r$ independent of $R,$

$(iii)$  if $(R_n)_{n \in \N}$ is a sequence of  real closed $(1,1)$-currents of order $0$ of uniformly bounded mass converging weakly to $R$ on $\B$ then $U_{R_n} \to U_R$ in $L^{1+1/(2k)}(\B_r).$  
\end{lemma}

\proof  This lemma is essentially classical. The new point is the estimate concerning $L^{1+1/(2k)}$-norm of the potential $U_R$ and its continuity in $R.$ These properties will be obtained by carefully examining steps in  the usual construction of $U_R,$ see \cite[p. 135]{Demailly_ag} for example.  

Let $R$ be a real closed $(1,1)$-current on $\B.$ Let $x\in \C^k$ be the canonical coordinate system. 
Let $\rho$ be a smooth function compactly supported in $\B$ and $\int_{\B}\rho dx=1.$   For $y \in \B,$ let $A_y: \B \to \B$ be the diffeomorphism defined by 
$$A_y(x):= x+\frac{1}{2}(1- \|x\|)y$$
for $x\in \B.$ Since $A_y$ is homotopic to $A_0:=\id$ through the homotopy $H_y: [0,1]\times \B \to \B$ defined by $H_y(t,x):=A_{ty}(x)$ for $t\in [0,1],$  the average 
$$R':=\int_{\B} (A_y^* R) \rho(y) dy$$
 is a smooth closed form  which is cohomologous to $R.$ Precisely, by the homotopy formula,  we have
$$R-R'= d L_1, \quad \text{where } L_1=L_1(R):= \int_B (H_y)_*([0,1]\otimes R)\rho(y) dy.$$ 
Observe that 
\begin{align}\label{ine-massSRphay}
\|R'\|_{L^\infty(\B)} \lesssim \|R\|_\B, \quad \|L_1\|_{\B}\lesssim \|R\|_\B.
\end{align}
Since $R'$ is a smooth closed form in $\B,$ we can use an explicit formula which we don't recall here (see \cite[p. 13]{Demailly_ag}) to  define a smooth form $L_2=L_2(R')$ on $\B$ such that 
$$R'= d L_2, \quad \|L_2\|_{L^\infty(\B)} \lesssim \|R'\|_{L^\infty(\B)}.$$
This combined with  (\ref{ine-massSRphay}) shows that for $L_3:=L_1+ L_2,$ we have
\begin{align}\label{ine-massSRphay2}
R= d  L_3, \quad \|L_3\|_\B \lesssim \|R\|_\B
\end{align}
and $L_3$ depends  continuously on  $R.$ Hence if $(R_n)_{n\in \N}$ is a sequence of $(1,1)$-current of order $0$ of uniformly bounded mass converging to $R$ then $L_3(R_n)$ is also of uniformly bounded mass and converges to $L_3(R).$  

Since $R$ is a real $(1,1)$-form,  $L_3$ is a real $1$-form.  We decompose $L_3$ into the sum of an $(1,0)$-form and  an $(0,1)$-form as 
\begin{align}\label{eqtachL3}
L_3= L_3^{(1,0)}+ L_3^{(0,1)}
\end{align}
 such that $L_3^{(1,0)}= \overline{ L_3^{(0,1)}}$ and $L_3^{(1,0)}, L_3^{(0,1)}$ are currents of order $0.$ We deduce from (\ref{ine-massSRphay2}) that 
\begin{align}\label{ine-massSRphay22}
 \|L_3^{(0,1)}\|_\B \lesssim \|R\|_\B
\end{align}
By a bidegree reason and the fact that $R= d L_3$, we get $\bar \partial  L_3^{(0,1)}=0.$ It is well-known that there is a distribution  $v$ define on an open neighborhood of $\overline \B_r$  with $\bar \partial v= L_3^{(0,1)}.$  We will recall briefly how to construct such a $v$ as a function of $L_3^{(0,1)}.$ The reference is \cite[p. 28]{Demailly_ag}.           
 
 Let $\rho$ be the function as above. We can assume $\rho\equiv 1$ on an open neighborhood of  $\overline \B_{r}.$  By the Koppelman formula, we have 
\begin{align}\label{eq-rhoS}
\rho L_3^{(0,1)}(x)= \bar \partial \int_\B K_1(x,y) \wedge \rho(y) L_3^{(0,1)}(y)+ \int_\B K_2(x,y)\wedge  \bar \partial \rho(y) \wedge L_3^{(0,1)}(y).
\end{align}
We don't recall the explicit formulae for $K_1,K_2$ but only emphasize that $K_1,K_2$ are the products of $\|x-y\|^{-2k+1}$ with smooth forms on $\C^k.$

Denote by $I_1, I_2$ the first and second integrals respectively of the right-hand side of (\ref{eq-rhoS}). We have 
$$\bar \partial I_1+ I_2= \rho L_3^{(0,1)}$$
which is equal to $L_3^{(0,1)}$ on $\B_r.$ 

By the singularity type of $K_1$ and the fact that $L_3^{(0,1)}$ is of order $0,$ we see that $I_1$ is a form with coefficients in $L^{1+1/(2k)}(\B)$ with 
\begin{align}\label{ine-L12kI1}
\|I_1\|_{L^{1+1/(2k)}(\B)}\lesssim  \|L_3^{(0,1)}\|_\B  \lesssim \|R\|_\B.
\end{align}
by (\ref{ine-massSRphay22}). On the other hand, since $\bar \partial \rho\equiv 0$ on an open neighborhood of  $\overline \B_r,$ the current $I_2$ is smooth on  $\B_{r'}$ for some $r'>r$.  Following exactly arguments  in \cite[p. 29]{Demailly_ag}, we obtain a smooth function $I_3$ on  $\B_{r'}$ for some $r'>r$ such that  $I_2= \bar \partial I_3$ on  $\B_r$ and 
 \begin{align}\label{ineI3chuanL12k}
\|I_3\|_{L^\infty(\B_r)} \le \|L^{(0,1)}_3\|_{\B} \lesssim \|R\|_\B
\end{align}
by (\ref{ine-massSRphay22})   and $I_3: R \mapsto I_3(R)\in L^\infty(\B_r)$ is continuous.   Thus if  $v:=(I_1+ I_3)$ then
$$L^{(0,1)}_3=\bar \partial v$$
on  $\B_r$. This  together with   (\ref{eqtachL3}) gives
$$L_3= \bar \partial v+ \partial \bar v.$$
We deduce from this and (\ref{ine-massSRphay2}) that 
$$R= d L_3= \partial \bar \partial(v - \bar v).$$
Hence $U_R:= 2\pi \Im v$ satisfies $R= \ddc U_R$  (recall $\ddc= (i/\pi) \partial \bar \partial$) and
$$\|U_R\|_{L^{1+1/(2k)}(\B_r)} \lesssim \|I_1\|_{L^{1+1/(2k)}(\B_r)}+ \|I_3\|_{L^{1+1/(2k)}(\B_r)}\lesssim \|R\|_\B$$
by (\ref{ine-L12kI1}) and (\ref{ineI3chuanL12k}). 

It remains to prove the continuity property of $U_R.$ We saw that $I_3, L_3$ are continuous in $R.$ We only need to check this property for $I_1.$  Let $(R_n)$ be the sequence as in the statement. We will show that  $I_1(R_n) \to I_1(R)$ in $L^{1+1/(2k)}(\B)$. By the above continuity property of $L_3,$ we have that  $S_n:=\rho L_3^{(0,1)}(R_n)$ is of uniformly bounded mass and converges to $S:=\rho L^{(0,1)}_3(R)$ as $n \to \infty.$   Write 
$$K_1(x,y)= \|x-y\|^{-2k+1} K'_1(x,y),$$
where $K'_1(x,y)$ is a smooth form.   For every small constant $\epsilon>0,$ let 
$$K_{1,\epsilon}(x,y):= \max\{\|x-y\|,\epsilon\}^{-2k+1} K'_1(x,y)$$
 which is a continuous form. As $\epsilon \to 0,$ we have $K_{1,\epsilon}(\cdot,y) \to  K_1(\cdot,y)$ in $L^{1+1/(2k)}(\B)$ uniformly in $y \in \B$. Thus  as $n \to \infty,$
$$\int_{\{y \in \B\}} \big(K_{1,\epsilon}(x,y) - K_1(x,y) \big) \wedge \big(S_n(y)- S(y)\big) \to 0$$
in $L^{1+1/(2k)}(\B)$ because the mass of $S_n$ is uniformly bounded. On the other hand, 
$$\int_{\{y \in \B\}} K_{1,\epsilon}(x,y) \wedge \big(S_n(y)- S(y)\big)$$ converges uniformly to $0$ as $\epsilon$ fixed because $K_{1,\epsilon}$ is continuous. We deduce that $I_1(R_n)\to I_1(R)$
in $L^{1+1/(2k)}(\B).$    The proof is finished. 
\endproof


Let $X$ be a  complex manifold.  A function from $X$ to $[-\infty, \infty)$ is said to be  \emph{quasi-p.s.h. function} if it can be written  locally the sum of a plurisubharmonic (p.s.h.) function and a smooth one. For every continuous $(1,1)$-form $\eta$, a quasi-p.s.h. function $\varphi$ is \emph{$\eta$-p.s.h.} if $\ddc \varphi+ \eta\ge 0.$ By partition of unity, every quasi-p.s.h. function is $\eta$-p.s.h. for some smooth form $\eta.$ For a given form $\eta,$ denote by $\PSH(\eta)$ the set of quasi-p.s.h. functions $\varphi$ for which $\ddc \varphi +\eta \ge 0.$ 

\begin{definition} A locally integrable function $\varphi$ on $X$ is said to be \emph{weakly d.s.h.} if $\ddc\varphi$ is a current of order $0$ on $X$.   Let $\cali{W}$ be the complex vector space of all weakly d.s.h. functions on $X.$ 
\end{definition}

Clearly, every quasi-p.s.h is weakly d.s.h..  A subset of $X$ is a  \emph{pluripolar set } if it is contained in $\{\varphi= -\infty\}$ for some quasi-p.s.h. function $\varphi$.  If $X$ is compact, every locally pluripolar set is pluripolar by  \cite{Vu_pluripolar}. In our proofs, we only  use a particular case of this result that every proper analytic subset of a compact manifold $X$ is pluripolar, see Lemma \ref{le_pluripola} below.

Consider now $X$ is \emph{compact.}  Let $\mu_0$ be a smooth probability measure on  $X.$  We use this measure to define $L^p$ norms on $X.$ For $\varphi \in \cali{W},$ put 
\begin{align}\label{norm_W}
\|\varphi\|_{\cali{W}}:= \bigg|\int_X \varphi d \mu_0\bigg|+ \|\ddc \varphi\|_X,
\end{align}
where $\|\cdot\|_X$ is the mass of a current on $X.$ We will write from now on $\|\cdot\|$ instead of $\|\cdot\|_X$ if no confusion arises. The function $\|\cdot\|_{\cali{W}}$ is a norm on $\cali{W}$ because if $\ddc \varphi=0$ then $\varphi$ must be a constant. The norm $\|\cdot\|_{\cali{W}}$ is similar to the norm of the space of d.s.h. functions on the K\"ahler case introduced by Dinh-Sibony \cite{DS_tm}. However, we don't know if these two norms are equivalent in this case.   

As in the K\"ahler case,   we introduce \emph{the topology} on $\cali{W}$ as follows: we say that $\varphi_n \in \cali{W}$ converges to $\varphi \in \cali{W}$  as $n \to \infty$ if $\varphi_n \to \varphi$ as currents and $\|\varphi_n\|_\cali{W}$ is uniformly bounded. 
We have the following crucial compactness result.

\begin{lemma} \label{le_compact}   Let $X$ be a compact complex manifold. There exists a constant $c$ so that for every weakly d.s.h. function $\varphi$ on $X$ with $\int_X \varphi d\mu_0=0,$ we have
\begin{align} \label{ine_L1ddcmass}
\|\varphi\|_{L^{1+1/(2k)}(X)} \le c \|\ddc \varphi\|_X.
\end{align}
Moreover,  given a positive constant $A,$  the set $\cali{W}_0$  of weakly quasi-p.s.h. functions $\varphi$ with $\int_X \varphi d \mu_0=0$  such that 
$\|\ddc \varphi\| \le A$
is compact in $L^{1+1/(2k)}(X).$  
\end{lemma}

A direct consequence of Lemma \ref{le_compact} is that if $\varphi_n \to \varphi$ in $\cali{W}$ then $\varphi_n \to \varphi$ in $L^{1+1/(2k)}.$  In K\"ahler case, a related version of the  inequality (\ref{ine_L1ddcmass}) for d.s.h. functions with $L^p$-norm  in place of $L^{1+1/(2k)}$-norm and $\|\cdot\|_*$  in place of $\|\cdot\|_X$ was proved in \cite{DS_tm} using cohomological tools for d.s.h. functions. Their proof uses cohomological arguments which are not applicable to prove   (\ref{ine_L1ddcmass}) for weakly quasi-p.s.h. functions. 

\proof    Consider a weakly quasi-p.s.h. function $\varphi$ with $\|\ddc \varphi\| \le A.$  Let $(W_j)$ be (finite) open covering of $X$ where $W_j$ are  local charts of $X$ biholomorphic to the unit ball of $\C^k.$   Since $\|\ddc \varphi\| \le A,$ by Lemma \ref{le-11current}, we have $\tau_j \in L^{1+1/(2k)}(W_j)$ for which $\ddc \tau_j= \ddc \varphi$ on $W_j$ and 
\begin{align}\label{inechuantauj}
\|\tau_j\|_{L^{1+1/(2k)}(W_j)} \lesssim A.
\end{align}
Hence, $\varphi -\tau_j$ can be represented by a pluriharmonic function on $W_j.$ For simplicity, we identify this function with $(\varphi -\tau_j).$ We deduce that $\varphi \in L^{1+1/(2k)}(X).$ 

We now suppose on the contrary that (\ref{ine_L1ddcmass}) does't hold, that means that there exists a sequence of  non-zero weakly quasi-p.s.h. functions $\varphi_n$ with $\int_X \varphi_n d\mu_0=0$ and 
$$\infty >\|\varphi_n\|_{L^{1+1/(2k)}(X)} \ge n \|\ddc \varphi_n\|_X.$$
By multiplying $\varphi_n$ by a positive constant,  we can assume that
\begin{align}\label{inechuantauj22}
\|\varphi_n\|_{L^{1+1/(2k)}(X)}=1.
\end{align}
 Thus, we get 
\begin{align} \label{ine_massTn}
\|\ddc \varphi_n\| \le 1/n.
\end{align}
Note that we still have $\int_X \varphi_n d\mu_0=0.$   Let $\tau_j^n$ be the function $\tau_j$ for $\varphi_n$ in place of $\varphi.$ Put $T_n:= \ddc \varphi_n.$ These currents of order $0$ are of uniformly bounded mass and converges to $0$ by (\ref{ine_massTn}). Lemma \ref{le-11current} tells us that $\tau^n_j$ converges to $0$ in $L^{1+1/(2k)}(W'_j),$ for every $W'_j \Subset W_j.$ We can also arrange that $(W'_j)$ is still a covering of $X.$ For simplicity, we can assume that $W'_j= W_j$ for every $j.$

Recall now that $\varphi_n - \tau^n_j$ is pluriharmonic on $W_j.$ The last function is of $L^{1+1/(2k)}$-norm bounded on $W_j$ because of (\ref{inechuantauj}) and (\ref{inechuantauj22}). The mean equality for pluriharmonic functions  implies that $(\varphi_n- \tau^n_j)$ is of $\cali{C}^l$-norm uniformly bounded on compact subsets of $W_j$  in $n \in \N$ for every $l \in \N.$ We deduce that by extracting a subsequence, we can suppose that $\varphi_n -\tau^n_j$ converging uniformly to a pluriharmonic function $\tau^\infty_j$ on  compact subsets of $W_j$ as $n \to \infty.$  Since $\|\tau^n_j\|_{L^{1+1/(2k)}(W_j)} \to 0,$ we obtain that  
$$\varphi_n \to \tau^\infty_j \quad \text{in } \, L^{1+1/(2k)}(W_j).$$
This yields that  function $\tau^\infty:= \tau^\infty_j$ on $W_j$ for every $j$ is a well-defined pluriharmonic function on $X.$ Since $X$ is compact, $\tau^\infty$ is a constant. This combined with $\int_X \varphi_n d\mu_0=0$ gives $\tau^\infty=0.$ We thus have proved that $\varphi_n \to 0$ in $L^{1+1/(2k)}(X),$ hence $\|\varphi_n\|_{L^{1+1/(2k)}}\to 0,$ a contradiction. Thus,  (\ref{ine_L1ddcmass}) holds. 

In order to prove the second desired assertion, we use again the function $\tau_j$ above. We have $\varphi -\tau_j$ is pluriharmonic on $W_j$ and by (\ref{ine_L1ddcmass}), the $L^{1+1/(2k)}$-norm of $\varphi$ is also $\lesssim A$. Thus, the $L^{1+1/(2k)}$-norm of  the pluriharmonic function $(\varphi -\tau_j)$ is $\lesssim A.$ It follows that its $\cali{C}^l$-norm is $\lesssim A$ as well. Hence, we can  extract a convergent subsequence of   $(\varphi -\tau_j)$ for $\varphi \in \cali{W}$ in $\cali{C}^l.$ This combined with the $L^{1+1/(2k)}$ continuity of $\tau_j$ in $T$ implies the desired assertion. The proof is finished.   
\endproof


We equip  the vector space $\cali{B}$ of Borel measurable functions on $X$ with the pointwise convergence topology: $h_n \to h$ if $h_n$ converges pointwise to $h$ almost everywhere (with respect to Lebesgue's measure).  Let $P$ be  a continuous linear endomorphism of the last vector space. Define $\cali{W}_P$ to be the set of $\varphi \in \cali{W}$ for which $P\varphi \in \cali{W}.$  

\begin{lemma}\label{le_compact2} There exists a constant $c$ such that 
\begin{align} \label{ine_comfcorWf}
\|P \varphi\|_{L^{1+1/(2k)}} \le c \big( \|\varphi\|_{\cali{W}}+ \|\ddc(P \varphi)\|\big),
\end{align}
for any $\varphi \in \cali{W}_P.$ In particular, there is a constant $c'$ so that 
\begin{align} \label{ine_comfcorWf2}
\|P \varphi\|_{L^{1+1/(2k)}} \le c \big(\|\ddc \varphi\|+ \|\ddc(P \varphi)\|\big)
\end{align}
for every $\varphi \in \cali{W}_P \cap \cali{W}_0.$ Moreover, if $\varphi_n \in \cali{W}_P \cap \cali{W}_0 \to \varphi$ as currents as $n\to \infty$ so that $\big(\|\ddc \varphi_n\|+ \|\ddc(P \varphi_n)\|\big)$ are uniformly bounded, then $P \varphi_n \to P \varphi$ in $L^{1+1/(2k)}.$  
\end{lemma}

\proof The inequality  (\ref{ine_comfcorWf2}) is a direct consequence of (\ref{ine_comfcorWf}) and Lemma \ref{le_compact}. 
Now suppose that there is a sequence $(\varphi_n) \subset \cali{W}_P$ for which 
\begin{align}\label{ine_contraryf*phi}
\|P \varphi_n\|_{L^{1+1/(2k)}}=1, \quad   \|\varphi\|_{\cali{W}}+ \|\ddc(P \varphi_n)\| \le 1/n.
\end{align}
Applying the compactness property in Lemma \ref{le_compact} to the sequence $(P \varphi_n)_{n \in \N}$, we see that by  extracting a subsequence of $\varphi_n$ if necessary, the sequence $P\varphi_n$ converges in $L^{1+1/(2k)}$ to a weakly d.s.h. function $\varphi'_\infty.$ Consequently,  
\begin{align} \label{eq_L1varphiinfry}
\|\varphi'_\infty\|_{L^{1+1/(2k)}}=1, \quad \|\ddc \varphi'_\infty\|=0.
\end{align}
Hence $\varphi'_\infty$ is a constant.   Since the convergence in $L^1$ implies the almost everywhere convergence of a subsequence, we can suppose also that $P \varphi_n$ converges almost everywhere to $\varphi'_\infty.$ 

On the other hand, the inequality of (\ref{ine_contraryf*phi}) allows us to use the compactness property in  Lemma \ref{le_compact} again for $(\varphi_n)$. Hence, we can extract a subsequence of $(\varphi_n)$ converging to $\varphi_\infty:=0$ in $L^{1+ 1/(2k)}$ and almost everywhere. Thus  $P\varphi_n$ converges almost everywhere to $P \varphi_\infty$ because of the continuity of $P.$   It follows that $\varphi'_\infty= P \varphi_\infty=0,$ note here $P(0)=0$ by the linearity of $P.$  This is a contraction because of (\ref{eq_L1varphiinfry}).  Thus (\ref{ine_comfcorWf}) follows. The last desired assertion follows directly from above arguments. 
The proof is finished.   
\endproof

Let $a \in \C^*$, $r$ a constant  in $(0, |a|)$ and $\delta>0$ a constant.  Assume that $P(1)=a,$ where $1$ is the constant function equal to $1$ on $X$.  Define $\cali{W}^\infty_{P,r,\delta}$ to be the set of all $\varphi \in \cali{B}$ such that  $P^n \varphi \in \cali{W}$ for every $n \ge 0$ and $$\|\ddc(P^n\varphi)\| \le  \delta  r^n$$
 for every $n \ge 0,$ here $P^0$ denote the identity map. By the linearity of $P,$ every constant function belongs to $\cali{W}^\infty_{P,r,\delta}.$   We equip $\cali{W}^\infty_{P,r,\delta}$ with the induced topology from that on $\cali{W}.$ Observe that $\cali{W}^\infty_{P,r,\delta}$ is closed in $\cali{W}$ and 
$$r^{-m} P^m(\cali{W}^\infty_{P,r,\delta}) \subset \cali{W}^\infty_{P,r,\delta}$$
 for every positive integer $m.$ Hence $\cali{W}^\infty_{P,r,\delta} \cap \cali{W}_0$ is compact and $P^m(\cali{W}^\infty_{P,r,\delta})$ is contained in the complex vector subspace $\tilde{\cali{W}}^\infty_{P,r,\delta}$ of $\cali{W}$ generated by $\cali{W}^\infty_{P,r,\delta}.$    

\begin{proposition}\label{pro_WPra}  
 There exists a continuous  linear functional $\mu_{P}: \tilde{\cali{W}}^\infty_{P,r,\delta} \to \C$ such that for every  complex measure $\nu$  with $L^{2k+1}$ density  on $X,$ $\nu(X)=1$ and    for every $\varphi \in \tilde{\cali{W}}_{P,r,\delta},$ we have 
\begin{align} \label{conver_numuP}
\big\langle  a^{-n} (P^n)_*  \nu,\varphi \big\rangle \to \langle \mu_{P},\varphi \rangle.
\end{align}
\end{proposition}

Here for $Q: \cali{B} \to \cali{B}$, by definition, $ \langle Q_* \nu, \varphi \rangle:= \langle \nu, Q \varphi \rangle$ for $\varphi \in \cali{B}$ such that $Q \varphi$ is $\nu$-integrable. 

\proof  Recall that  $\mu_0$  is a smooth probability volume form on $X.$ We only need to  construct $\mu_P$ on $\cali{W}^\infty_{P,r,\delta}$ and prove (\ref{conver_numuP}) for $\varphi \in \cali{W}^\infty_{P,r,\delta}.$ The extension of $\mu_P$ to $\tilde{\cali{W}}^\infty_{P,r,\delta}$ is automatically done by using the linearity of $(P^n)_*  \nu$ and (\ref{conver_numuP}). 

 Let $\varphi \in \cali{W}^\infty_{P,r,\delta}.$  Put  $b_0:= \int_X \varphi d \mu_0$ and $\varphi_0:= \varphi -b_0.$  We define two sequences $\varphi_n, b_n$ as follows. Put
$$b_n= b_n(\varphi):= \int_X (P\varphi_{n-1}) d\mu_0,  \quad \varphi_n:= P\varphi_{n-1}- b_{n}$$
 for $n\ge 1.$ We have  $r^{-n} \varphi_n \in \cali{W}_0 \cap \cali{W}^\infty_{P,r,\delta}$ and $\ddc(P^m \varphi_n)= \ddc (P^{m+n} \varphi)$ for every $n,m.$     By  Lemma \ref{le_compact2}, we have 
$$ \|\varphi_n\|_{L^{1+1/(2k)}} \le c\big((\|\ddc (P\varphi_{n-1})\|+ \|\ddc \varphi_{n-1}\|\big),\quad |b_n|  \le c \big(\|\ddc (P\varphi_{n-1})\|+ \|\ddc \varphi_{n-1}\|\big).$$
for some constant $c$ independent of $n,\varphi.$ It follows that 
\begin{align} \label{ine_bnPP}
\|\varphi_n\|_{L^{1+1/(2k)}} \le  c \big(\|\ddc (P^n\varphi)\|+  \| \ddc (P^{n-1}\varphi) \|) \le c \delta (r+1) r^{n-1}, \quad |b_n|  \le c \delta (r+1) r^{n-1}
\end{align}
for $n\ge 1.$  Since $P(1)=a$ we have $P(b_n)= a b_n$ for every $n.$  Using this gives
\begin{align}\label{eq_vietfnsaophiPP}
a^{-n} P^n \varphi &= b_0+ a^{-n}P^n \varphi_0= b_0+  a^{-n}P^{n-1}( P \varphi_0)\\
\nonumber
&= b_0+ a^{-1}b_1+ a^{-n}P^{n-1} \varphi_1=\\
\nonumber
&= \cdots = b_0+ a^{-1}b_1+ \cdots a^{-n} b_{n}+ a^{-n}\varphi_n. 
\end{align}
Put $b'_n= b'_n(\varphi):=b_0+ a^{-1}b_1+ \cdots a^{-n} b_{n}$  which converges to a number $b'_\infty$ (depending on $\varphi$)  by (\ref{ine_bnPP}) and the fact that $|a|>r$.  We deduce from (\ref{eq_vietfnsaophiPP}) that 
$$|a^{-n}P^n \varphi -b'_n | \le  |a|^{-n} |\varphi_n|.$$
This combined with the first inequality of (\ref{ine_bnPP}) implies that $a^{-n}P\varphi$ converges to $b'_\infty$ in $L^{1+1/(2k)}.$ Precisely, we have 
\begin{align}\label{ine_L1jbnphayPP}
\|a^{-n}P^n \varphi -b'_n \|_{L^{1+1/(2k)}} \lesssim \delta  |a|^{-n} r^{n}.
\end{align}
Since $\nu(X)=1,$ we get 
$$\langle  a^{-n} (P^n)_* \nu, \varphi \rangle - b'_n= \langle  \nu,  a^{-n} P^n\varphi - b'_n \rangle.$$
Using this, (\ref{ine_L1jbnphayPP}) and  H\"older's inequality implies  that
$\langle  a^{-n} (P^n)_* \nu, \varphi \rangle$ converges to $b'_\infty=b'_\infty(\varphi)$ because $\nu$ has $L^{2k+1}$ density.  Define $\langle \mu_P, \varphi\rangle := b'_\infty(\varphi) $ which is independent of $\nu.$ We then obtain  the desired convergence toward $\mu_P.$ 

Consider a sequence $\tilde{\varphi}_m \to \varphi$ in $\cali{W}^\infty_{P,r,\delta}.$ Let $\tilde{b}_{n m}, \tilde{\varphi}_{n m}$ be respectively the $b_n$ and $\varphi_n$   for $\tilde{\varphi}_m$ in place of $\varphi.$   By the last assertion of    Lemma \ref{le_compact2}, $\tilde{b}_{n m} \to b_n $ as $m \to \infty$ for every $n$ and  (\ref{ine_bnPP}) still holds for $\tilde{b}_{n m}, \tilde{\varphi}_{n m}$ in place of $b_n,\varphi_n$.  We infer that $\tilde{b}'_{n m} \to b'_n$ and $a^{-n} \tilde{\varphi}_{n m} \to 0$ in $L^{1+1/(2k)}$ as $m \to \infty.$ Thus, $\langle \mu_P, \tilde{\varphi}_m \rangle \to \langle \mu_P, \varphi\rangle$ as $m \to \infty.$  In other words, $\mu_P$ is continuous. 
The proof is finished. 
\endproof

Let $X$ be a compact complex manifold and $f$ a meromorphic self-map on $X.$  Denote by $\Gamma$ the graph of $f$ on $X \times X$ and $\pi_1,\pi_2$ the restrictions to $\Gamma$ of the natural projections from $X\times X$ to the first and second components respectively.    

Let $\Phi$ be a form with measurable coefficients on $X.$ We say that $\Phi \in L^1$ if its coefficients are $L^1$ functions (with respect to Lebesgue's measure on $X$).   If  $\Omega$ is an open  Zariski dense subset of $X$ such that $\pi_2$ is an unramified covering over $\Omega,$  the form $f_* \Phi:= (\pi_2|_{\pi_2^{-1}(\Omega)})_*(\pi^*_1 \Phi)$ is a  measurable form on $\Omega.$ Hence  $f_*\Phi$ is a measurable form on $X$ independent of $\Omega.$ We can check that    $f_*: \cali{B} \to \cali{B}$ is continuous. Consequently, $f_*$ is an example of the map $P$ considered above. 

If  $f_* \Phi \in L^1,$ then   we can define $f_* \Phi$ to be a current of order $0$ induced by $f_* \Phi$ on $X$.  This definition is independent of the choice of $\Omega.$  Note that the pull-back by $f$ of smooth functions or smooth forms is always in $L^1$.  The following is similar to results  in  \cite{Meo,DinhSibony_pullback}.

\begin{lemma} \label{le_quasipull} For every quasi-p.s.h. function $\varphi$ on $X,$ we have $f_* \varphi \in L^1$ and if $\ddc\varphi + \eta \ge 0$ for some continuous $(1,1)$-form $\eta>0,$ then $\ddc (f_* \varphi)+ f_* \eta \ge 0.$ In particular, 
\begin{align} \label{inclu_fsao}
(f^n)_* \varphi \in \cali{W}_{f_*} \cap \cali{W}.
\end{align}
\end{lemma}  
\noindent    
The inclusion (\ref{inclu_fsao}) explains the crucial roles of $\cali{W}_{f_*}, \cali{W}$ in our study.

\proof Let $\sigma:\Gamma' \to \Gamma$ be a desingularisation of $\Gamma.$  Let $\Omega$ be as above. Put $\pi_j':= \pi_j \circ \sigma$ for $j=1,2.$ Since $\varphi$ is quasi-p.s.h.,  $\varphi \circ \pi_1'$ is so. Thus, $\varphi \circ \pi_1= \sigma_*(\varphi \circ \pi_1')$ is in $L^1(\Gamma_f).$ Since  
$$\|f_* \varphi\|_{L^1(\Omega)}= \|(\pi_2)_* (\varphi \circ \pi_1)\|_{L^1(\Omega)} \lesssim\|\varphi \circ \pi_1\|_{L^1(\Gamma)},$$ 
we get the first desired assertion.

By \cite{Blocki_kolo_regular} and the fact that $\eta>0$, there exists  a decreasing sequence of smooth quasi-p.s.h functions $\varphi_n$ converging pointwise to $\varphi$ such that $\ddc \varphi_n + \eta \ge 0$ for every $n$.   By Lebesgue's dominated convergence theorem, the sequence $\varphi_n \circ \pi_1'$ converges in $L^1$ to $\varphi \circ \pi_1'.$ It follows that the sequence of positive smooth forms $\ddc (\varphi_n \circ \pi_1')+ \pi_1'^* \eta$  converges weakly to $\ddc (\varphi \circ \pi_1')+ \pi_1'^* \eta.$ Thus, the last current is also positive.  Now observe that  
$$(\pi_2')_*(\ddc (\varphi \circ \pi_1')+ \pi_1'^* \eta)=\ddc \big((\pi_2')_* \pi_1'^* \varphi \big)+ (\pi'_2)_*\pi_1'^* \eta=\ddc \big((\pi_2)_*\pi_1^*\varphi\big)+ (\pi_2)_*\pi_1^* \eta$$
because $\pi_1^* \varphi$ and $\pi_1^* \eta$ have no mass on sets of Lebesgue measure zero. Thus $\ddc (f_* \varphi)+ f_* \eta \ge 0.$ 

Note that $f_* \eta$ has  finite mass on $X.$ We infer that $f_* \varphi \in \cali{W}.$ In other words, $\varphi \in \cali{W}_{f_*} \cap \cali{W}.$ Applying this  to $f^n$ instead of $f$ and using the formula  that $(f^n)_* \varphi= f_* (f^{n-1})_* \varphi$ as functions  on some suitable open dense subset of $X,$ we obtain (\ref{inclu_fsao}).  This finishes the proof.
\endproof


\begin{lemma} \label{le_pshpush} Let $X$ be a compact complex manifold of dimension $k$  and $f$ a meromorphic self-map on $X.$ Let $\varphi$ be a quasi-p.s.h. function on $X$ with  $\ddc \varphi+ \eta \ge 0$ for some continuous $(1,1)$-form $\eta.$ Then given every positive constant $\epsilon,$ there exists a constant $c_\epsilon$ independent of $\varphi,\eta$ for which
\begin{align} \label{ine_massfnsapvarphi}
\|\ddc (f^n)_*\varphi\| \le c_\epsilon (d_{k-1}(f)+\epsilon)^n \|\eta\|_{L^\infty}
\end{align}
for every $n \ge 1.$ 
\end{lemma}

\proof  By replacing $\eta$ by a strictly positive smooth form dominating it, we can assume that $\eta>0.$   Let $\omega$ be a Gauduchon metric on $X,$ that means that $\omega$ is a Hermitian metric and  $\ddc \omega^{k-1}=0,$ see \cite{Gauduchon}.  Let $\Gamma_n$ be the graph of $f^n$ and $\pi_{1,n}, \pi_{2,n}$ the natural maps from $\Gamma_n$ to the first and second components of $X\times X.$  
By Lemma \ref{le_quasipull}, the current $\ddc (f^n)_*\varphi + (f^n)_* \eta$ is positive. Thus, using $\ddc \omega^{k-1}=0$ gives  
$$\| \ddc (f^n)_*\varphi + (f^n)_* \eta \| \lesssim \langle \ddc (f^n)_*\varphi + (f^n)_* \eta, \omega^{k-1} \rangle =  \langle (f^n)_*\eta, \omega^{k-1} \rangle \lesssim \langle (f^n)_*\omega, \omega^{k-1} \rangle$$
This combined with the  definition of $d_{k-1}(f)$ gives 
\begin{align*}
\|\ddc (f^n)_*\varphi+ (f^n)_* \eta\| \le c_\epsilon (d_{k-1}(f)+\epsilon)^n \|\eta\|_{L^\infty}.
\end{align*}
The desired inequality then follows immediately.  The proof is finished. 
\endproof


We now come to the end of the proof of the first main result. 

\begin{proof}[End of Proof of Theorem \ref{the_topodereedi}]  Fix a positive constant $\epsilon$ for which $d_k > d_{k-1}+\epsilon.$   Put 
$$P:= f_*, \quad  a:= d_k, \quad  r:= (d_{k-1}+\epsilon), \quad \delta:=c_\epsilon,$$
 where $c_\epsilon$ is the constant in Lemma \ref{le_pshpush}.  Let $\varphi$ be a quasi-p.s.h. with $\ddc \varphi+ \eta \ge 0$ for some continuous $(1,1)$-form $\eta>0$ so that  $\|\eta\|_{L^\infty}\le 1.$ We have $P(1)=a$ and $\varphi \in \cali{W}^\infty_{P,r,\delta}$  by Lemma \ref{le_pshpush}. Every quasi-p.s.h. function is in $\tilde{\cali{W}}^\infty_{P,r,\delta}.$     Since $\nu$ has no mass on proper analytic subsets of $X,$  observe that 
\begin{align} \label{eq_chuyenpushfvarphhi}
\langle (f^n)^*\nu, \varphi \rangle =\langle  \nu, (f^n)_*\varphi \rangle= \langle  \nu,  P^n \varphi \rangle
\end{align}
because we only need to consider integrals on an open Zariski dense subset of $X.$ Applying Proposition \ref{pro_WPra} to $P$, we obtain a continuous functional $\mu_P$ on $\tilde{\cali{W}}^\infty_{P,r,\delta}$ such that 
$$\langle d_k^{-n}(f^n)^*\nu, \varphi \rangle \to \langle \mu_P, \varphi \rangle,$$
for every $\varphi \in \tilde{\cali{W}}^\infty_{P,r,\delta}.$  By choosing $\nu \ge 0,$ we see that  $\langle \mu_P, \varphi \rangle \ge 0$ if $\varphi \ge 0.$   Let $\mu_f$ be the  probability measure on $X$ defined by $\langle \mu_f, \varphi \rangle:=\langle \mu_P, \varphi \rangle$ for every smooth function $\varphi$. Recall here that smooth functions are quasi-p.s.h. on $X.$    We will prove that $\mu_f= \mu_P$ for every quasi-p.s.h. function $\varphi.$

Consider a sequence of smooth  quasi-p.s.h. functions $\varphi'_n$ with $\ddc \varphi'_n +\eta\ge 0$ decreasing to $\varphi,$ we get  $\langle \mu_f, \varphi'_n \rangle =\langle \mu_P, \varphi'_n \rangle  $ and $\langle \mu_f, \varphi'_n \rangle \to \langle \mu_f, \varphi \rangle$ by Lebesgue's monotone convergence theorem. This combined with the continuity of $\mu_P$ gives  $\langle \mu_f, \varphi \rangle =\langle \mu_P, \varphi \rangle$. Thus  we get 
\begin{align} \label{conver_psh}
\lim_{n \to \infty} \langle d_k^{-n}(f^n)^* \nu - \mu_f, \varphi \rangle =0
\end{align}
for every quasi-p.s.h. function $\varphi$ on $X.$   

Since quasi-p.s.h functions  are $\mu_f$-integrable,  $\mu_f$ has no mass on pluripolar sets. By Lemma \ref{le_pluripola}  below,  proper analytic subsets of $X$ are  pluripolar.   This implies that $\mu_f$ has no mass on proper analytic subsets of $X.$ We deduce  that  the pull-back $f^* \mu_f$ is well-defined. Here we only take the  pull-back of  $\mu_f$ on an open Zariski subset $\Omega$ of $X$ where $\pi_2$ is an unramified covering. One can check that this definition is independent of the choice of $\Omega$ and if $(\Phi_n)_{n \in \N}$ is  a sequence of positive measures having no mass on proper analytic subsets of $X$ and converging to $\mu_f$, then $f^* \Phi_n$ converges to $f^* \mu_f$ because the mass of $f^* \Phi_n$ converges to that of $f^* \mu_f,$  see for example \cite[Le. 3.6]{Vu_intersec}.    The equality 
\begin{align}\label{eq_totallyinvariantmu}
d^{-1}_k f^* \mu_f =\mu_f
\end{align}
is obtained by applying the pull-back $f^*$ to the convergence $d_k^{-n}(f^n)^* \nu \to \mu_f,$ where $\nu$ is a smooth probability measure.   Since we have $f_* f^*=d_k$ on Borel measurable functions, we get $f_* \mu_f = \mu_f,$ in other words, $\mu_f$ is invariant by $f$. 

Let $I_f$ be the indeterminancy set of $f.$  Put  $Z: =\cup_{n \in \Z} f^n(I_f).$ The measure $\mu_f$ has no mass on $Z.$  The entropy of $\mu_f$ is by definition that of $\mathbf{1}_{X \backslash Z}\mu_f$ with respect to $f|_{X \backslash Z}.$    By an inequality of Parry \cite{Parry, DS_book}, using $f^*\mu_f=d_k \mu_f,$ we deduce that  the entropy of $\mu_f$ is at least $\log d_k.$  

Assume now  $f$ is holomorphic.  To prove that $\mu_f$ is H\"older continuous on $\PSH(\omega),$ we use a known idea from \cite{DS_book}.  Without loss of generality, we can assume that $\|\omega\|_{L^\infty} \le 1.$   Let $\varphi, \psi$ be two quasi-p.s.h. functions in $\PSH(\omega)$. Recall that they are in   $\cali{W}^\infty_{P,r,\delta}.$ 

Let $b_n(\varphi), b_n(\psi)$ be as  in the proof of Proposition \ref{pro_WPra}.  Let $J_f$ be the Jacobian of $f.$  We have 
$$\|f_* \varphi - f_* \psi\|_{L^1}= \sup_{\|h\|_{L^\infty} \le 1} | \langle f_* \varphi - f_* \psi, h \mu_0 \rangle|=\sup_{\|h\|_{L^\infty} \le 1} | \langle \varphi - \psi, (h\circ f) f^*\mu_0 \rangle|$$
which is  
$$\le  \| J_f\|_{L^\infty} \|\varphi - \psi\|_{L^1}.$$
Applying the last inequality for $f^n$ in place of $f$ gives
$$|b_n(\varphi)- b_n(\psi)| \le  2^n \| J_f\|^n_{L^\infty} \|\varphi - \psi\|_{L^1}.$$ 
 Put
$$A_1:=  \sum_{n=0}^{M+1} d_k^{-n}[b_n(\varphi)- b_n(\psi)], \quad A_2:=  \sum_{n=M+1}^\infty d_k^{-n} [b_n(\varphi)- b_n(\psi)].$$
Using (\ref{eq_vietfnsaophiPP}) gives 
$$\langle \mu_f, \varphi - \psi\rangle = A_1+ A_2,  \quad |A_1| \le  \sum_{n=0}^M d_k^{-n} 2^n \|J_f\|^n_{L^\infty} \|\varphi -\psi\|_{L^1}, \quad |A_2| \lesssim  (d_{k-1}+\epsilon)^M d_k^{-M}.$$ 
Consider the case where $2 \|J_f\|_{L^\infty} \le d_k.$ We have $|A_1| \le M \|\varphi -\psi\|_{L^1}.$ By choosing $M$ to be smallest integer for which $M \ge -\log \|\varphi - \psi\|_{L^1}/ \log \tau,$ where $\tau:= d_k/(d_{k-1}+ \epsilon),$ for every constant $\epsilon>0,$ we obtain that 
$$|\langle \mu_f, \varphi - \psi\rangle| \le  |A_1|+| A_2| \lesssim \|\varphi -\psi\|_{L^1}^{1-\epsilon}$$
which implies that $\mu_f$ is H\"older continuous in this case. It remains to treat the case  $2\|J_f\|_{L^\infty} \ge d_k.$ We have 
$$|A_1| \le M 2^M d_k^{-M} \|J_f\|_{L^\infty}^M \|\varphi -\psi\|_{L^1}+ \tau^{-M}.$$
Choose $M:= - \log \|\varphi -\psi\|_{L^1}/ \log (2 d^{-1}_k \tau \|J_f\|_{L^\infty}).$ We see that 
$$|A_1|+ |A_2| \lesssim - \log \|\varphi -\psi\|_{L^1}  \|\varphi -\psi\|_{L^1}^{\log \tau/ \log (2d^{-1}_k \tau \|J_f\|_{L^\infty})}.$$
Hence, $\mu_f$ is also H\"older continuous in this case.  This finishes the proof.
\end{proof}

Now we would like to say some words about Theorem \ref{the_mixing}. If one tries to mimic the arguments in the proof of \cite[The. 1.3]{DS_tm} to prove Theorem \ref{the_mixing}, we are led to estimating $|\langle \mu_f, |\varphi_n| \rangle |.$ The measure $\mu_f$ still satisfies the property that for every $\omega$-p.s.h. function $\varphi$ with $\sup_X \varphi=0$ is of $L^1(\mu_f)$-norm uniformly bounded, see \cite[Pro. 2.3]{DS_tm}. But unlike the K\"ahler case, we don't know whether $\varphi_n$ is the difference of two $\omega$-p.s.h functions. So this explains why we cannot apply directly the approach in \cite{DS_tm} to get a decay of correlation for $\mu_f.$ 
 


\begin{lemma} \label{le_pluripola} Every proper analytic subset $V$ of a compact complex manifold $X$ is a pluripolar set  on $X$.
\end{lemma}

\proof We use here the idea in \cite{DS_tm} where the authors prove the same result when $X$ is K\"ahler.  Suppose now that $V$ is smooth and $\codim V \ge 2$ (since otherwise the problem is trivial). Let $\sigma: \widehat X\to X$ be the blowup of $X$ along $V.$ Denote by $\widehat V$ the exceptional hypersurface.

Let $\omega$ be a positive definite Hermitian form on $X.$  Let $\widehat\omega_h$ be a Chern form of $\mathcal{O}(-\widehat V)$ whose restriction to each fiber of  $\widehat{V} \approx \P(E)$ is strictly positive.    By scaling $\omega$ if necessary, we can assume that  $\widehat \omega:=  \sigma^* \omega+ \widehat \omega_h>0.$   Since $\sigma_* \widehat \omega_h= \sigma_* \widehat \omega-  \omega,$ the closed current $\sigma_* \widehat \omega_h$ is quasi-positive. Thus there exists a quasi-p.s.h. function $\varphi$ on $\widehat X$ such that  
\begin{align}\label{eq_bieuthucomegah}
\sigma_* \widehat \omega_h= \ddc \varphi + \eta
\end{align}
 for some smooth closed form $\eta.$  By multiplying $\widehat \omega_h$ by a strictly positive constant, we have  $\sigma^* \sigma_* \widehat \omega_h= \widehat \omega_h+  [\widehat V].$   Thus   $\big| \varphi \circ \sigma(\widehat x) - \log \dist(\widehat x, \widehat V)|$  is a bounded function on  $\widehat X$. As a consequence,  
\begin{align} \label{ine_singvarphi}
|\varphi(x) - \log\dist(x, V)|  \lesssim 1
\end{align}
on compact subsets of $X.$ Consequently, $V$ is contained in $\{\varphi= -\infty\}.$  Hence $V$ is pluripolar in this case. 

By the above construction,  we can construct  a Hermitan metric on the blowup $\widehat X$ of $X$ along  $V$ as the sum of a pull-back of a Hermitian one on $X$ and a suitable Chern form of $\mathcal{O}(-\widehat V).$ Hence, if  $\sigma':\widehat X' \to X$ is a composition of  blowups along smooth submanifolds, then there are a smooth closed $(1,1)$-form $\eta'$ on $\widehat X'$ and a Hermitian metric $\omega$ on $X$ such that  $\widehat \omega'= \sigma'^* \omega+ \eta'$ is a Hermitian metric on $\widehat X'$. 

Consider now the general situation where $V$ is an analytic subset of $X$. Since a finite union of pluripolar sets is again pluripolar,   it is enough to prove that the regular part $\Reg V$ of $V$ is a pluripolar set because we can write $V$ as a finite union of the regular parts of suitable analytic subsets of $X.$   By Hironaka's desingularisation, there is a composition $\sigma': \widehat X' \to X$ of blowups along smooth submanifolds  which don't intersect $\Reg V$ (or its inverse images) such that the strict transform $\widehat V'$ of $V$ is  smooth. 

Let $\widehat \omega', \omega, \eta$ be as above.  By the above arguments,  $\widehat V'\subset \{\widehat \varphi'= -\infty\}$ for some  quasi-p.s.h. function $\widehat \varphi'$ on $\widehat X'$  and $\ddc \widehat \varphi' + \widehat \omega' \ge 0.$  Put $S:= \sigma'_*(\ddc  \widehat \varphi'+ \eta')$ which is a closed $(1,1)$-current  on $X$ and $S + \omega\ge 0.$ We can write $$S= \ddc \varphi_S+ \eta_S, \quad \sigma_* \eta'= \ddc \psi+ \eta$$
for some smooth closed forms $\eta_S, \eta.$ We have 
$$\ddc \varphi_S+ \eta_S+ \omega \ge 0, \quad  \ddc \psi+ \eta+ \omega \ge 0.$$
 Thus $\varphi_S, \psi$ are  quasi-p.s.h. functions on $X.$ Moreover,  we also have
$$\varphi_S= \sigma'_*(\widehat  \varphi')+ \psi+ \text{ a smooth function}$$
on an open neighborhood of $\Reg V$ on which $\sigma'$ is biholomorphic.  Consequently, $\Reg V \subset \{\varphi'_S= -\infty\}.$  This finishes the proof. 
\endproof

\section{The set $W^{1,2}_{*,f}$} \label{secW12}

In this section, we prove Theorem \ref{the_mixing}. 
Our idea is to consider suitable test functions in the Sobolev space $W^{1,2}.$ This approach is inspired by \cite{DS_decay}.  

Fix a smooth volume form $\mu_0$ on $X$ and we use this form to  define the norm on the  space $L^2(X).$    Let $W^{1,2}$ be the space of real-valued function $\varphi \in L^2(X)$ such that $d \varphi$ has $L^2$ coefficients. Recall the following Poincar\'e-Sobolev inequality: for $\varphi \in W^{1,2}$ with $\int_X \varphi d\mu_0 =0,$ we have 
\begin{align} \label{ineSP}
\| \varphi \|_{L^2} \le c \| d \varphi\|_{L^2}, 
\end{align}
for some constant $c$ independent of $\varphi,$ see for example \cite[Pro. 3.9]{Hebey_sobolev} or  \cite{Evans_pde}.  
Observe that  the term $\|d \varphi\|^2_{L^2}$ is comparable with the mass of the positive current $i \partial \varphi \wedge \bar \partial \varphi.$   We have the following lemma. 

\begin{lemma} \label{le-tinhchuanW12} (\cite[Pro. 3.1]{DS_decay})  Let $I$ be a compact subset of $X$ of the Hausdorff $(2k-1)$-dimensional measure zero. Let $\varphi$ be a real-valued  function  $L^1_{loc}(X \backslash I).$ Assume that the coefficients of $d \varphi$  are in $L^2(X \backslash I).$ Then $\varphi \in W^{1,2}$ and there exists a compact  subset $M$ of $X \backslash I$ and a constant $c>0$  both independent of $\varphi$ such that 
$$\| \varphi \|_{L^1(X)} \le c (\|\varphi\|_{L^1(M)}+\| d \varphi\|_{L^1(X)}).$$ 
\end{lemma}


    
The following is the central object in this section. 

\begin{definition} Let $W^{1,2}_{*,f}$ be the subset of $W^{1,2}$ consisting of $\varphi$ such that there exist $m_1\in \N,$ a continuous $(1,1)$-form $\eta$ and an $\eta$-p.s.h. function $\psi$  satisfying 
\begin{align}\label{ine-defW12f}
i \partial \varphi \wedge \bar \partial \varphi \le \ddc \big((f^{m_1})_*\psi\big)+ (f^{m_1})_*\eta
\end{align}
as currents.  \emph{A size representative} of $\varphi$ is $\mathbf{m}:=(m_0, m_1)$, where $m_0$ is an upper bound of  $\|\eta\|_{L^\infty}$. 
\end{definition}
If $X$ is K\"ahler, $W^{1,2}_{*,f}$ coincides with  the space  $W^{1,2}_*$  considered in \cite{DS_decay} which is independent of $f$. In that context, the space $W^{1,2}_*$ is studied in details in \cite{Vigny} and used in \cite{DLW} for the study of correspondences on Riemann surfaces with two equal dynamical degrees.    Let $\epsilon$  be a strictly positive constant such that $d_{k-1}+\epsilon < d_k$.  We have the following observation.

\begin{lemma} \label{le-sizereprese} Let $\varphi \in W^{1,2}_{*,f}$ and $\mathbf{m}=(m_0,m_1)$ a size representative of $\varphi.$ Then we have 
$$\|d \varphi\|_{L^2} \le c_\epsilon m_0^{1/2}  (d_{k-1}+\epsilon)^{m_1/2}$$
for some constant $c_\epsilon$ independent of $\varphi.$
\end{lemma} 

\proof Let $\eta$ be as in  (\ref{ine-defW12f}).  Let $\omega$ be a Hermitian metric on $X$ with  $\ddc \omega^{k-1}=0$. By testing  $ \ddc \big((f^{m_1})_*\psi\big)+ (f^{m_1})_*\eta$ with this form, we see that the norm of  $ \ddc \big((f^{m_1})_*\psi\big)+ (f^{m_1})_*\eta$ is equal to $\int_X (f^{m_1})_* \eta \wedge \omega^{k-1}$
which is bounded by  $c_\epsilon m_0 (d_{k-1}+\epsilon)^{m_1}$ for some constant $c_\epsilon$ independent of $\eta, m_0, m_1$.    The desired inequality then follows. This finishes the proof. 
\endproof

Let $\varphi \in W^{1,2}_{*,f}.$ Define $\varphi^+:= \max\{\varphi, 0\}$ an $\varphi^-:= \max\{-\varphi, 0\}.$ Consider a Lipschitz function $\chi: \R \to \R.$ We have
$\partial (\chi \circ \varphi)= (\chi' \circ \varphi) \partial \varphi.$ This can be seen by using a sequence of smooth functions converging to $\varphi$ in $W^{1,2}.$ We deduce that 
$$i \partial (\chi \circ \varphi) \wedge  \bar\partial (\chi \circ \varphi)= (\chi' \circ \varphi)^2 i\partial \varphi \wedge  \bar\partial \varphi.$$
Consequently, $\chi \circ \varphi \in W^{1,2}_{*,f}.$ In particular, by letting $\chi(t):= |t|$,  $\max\{t, 0\}$ or $\max\{-t,0\}$ for $t \in \R,$ we obtain the following crucial property.  

\begin{lemma} \label{le_repretrituyedoi} For every $\varphi \in W^{1,2}_{*,f},$ if  $\mathbf{m}=(m_0,m_1)$ is a size representative of $\varphi,$ then $\mathbf{m}$ is also a size representative of $|\varphi|,$ $\varphi^+$ and $\varphi^-.$ 
\end{lemma}

We already know that the pushforward of a quasi-p.s.h. function by $f$ is a weakly d.s.h. function. The following result, which explains the role of  $W^{1,2}_*$ in our study, gives a more precise description in the case of bounded quasi-p.s.h. functions.

\begin{lemma} \label{leppushvounded}  Every bounded quasi-p.s.h. function is in $W^{1,2}_{*,f}$ and $f_*$ preserves $W^{1,2}_{*,f}.$ Moreover, for every $\varphi \in W^{1,2}_{*,f},$ if  $\mathbf{m}=(m_0,m_1)$ is a size representative of $\varphi,$ then $\mathbf{m}':= (d_k m_0, m_1+1)$ is a size representative of $f_* \varphi$ and 
\begin{align}\label{inefsaofW12}
\|f_* \varphi\|_{L^2} \le c(\|\varphi\|_{L^1}+ \|d(f_*\varphi)\|_{L^2})
\end{align}
for some constant $c$ independent of $\varphi. $
\end{lemma}

\proof  Let $\varphi$ be a bounded quasi-p.s.h. function and $f: X \to X$ a dominant meromorphic map.  Using the identity 
$$2i \partial \varphi \wedge \bar \partial \varphi=i\partial \bar \partial \varphi^2- 2 \varphi i \partial \bar \partial \varphi$$
we see that there exist a continuous $(1,1)$-form $\eta$ and  an $\eta$-p.s.h. function $\psi$ for which $i \partial \varphi \wedge \bar \partial \varphi \le \ddc \psi+ \eta.$ Hence $\varphi \in W^{1,2}_{*,f}.$

Now let $\varphi$ be an arbitrary element of $W^{1,2}_{*,f}.$ Let $\eta$ and  $\psi$ be such that (\ref{ine-defW12f}) holds.    Fix an open Zariski dense subset $\Omega$ of $X$ on which $f_* \varphi, (f^{m_1})_* \psi, (f^{m_1})_* \eta$ are  well-defined functions or forms and $\pi_1$ is a unramified covering on $f^{-1}(\Omega)$. We have $f_* \varphi \in L^1_{loc}(\Omega)$ and 
\begin{align}\label{inechuanL1fsaoW12}
\|f_* \varphi\|_{L^1{(K)}} \le c \|\varphi\|_{L^1},
\end{align}
for any compact $K$ in $\Omega$ and  some constant $c$ independent of $\varphi.$    Note that $X \backslash \Omega$ is a proper analytic subset of $X,$ hence, is of Hausdorff $(2k-1)$-dimensional measure zero. On $\Omega,$ by the Cauchy-Schwarz inequality, we have
\begin{align*}
i \partial(f_*\varphi) \wedge \bar \partial (f_*\varphi)&\le d_k f_*(i \partial \varphi \wedge \bar \partial \varphi) \le d_k f_*\big[\ddc \big((f^{m_1})_*\psi\big)+ (f^{m_1})_*\eta\big]\\
&=d_k[\ddc \big((f^{m_1+1})_*\psi\big)+ (f^{m_1+1})_*\eta].
\end{align*}
It follows that $d (f_*\varphi) \in L^2(\Omega).$ By this and Lemma \ref{le-tinhchuanW12}, we get $f_* \varphi \in W^{1,2}.$ Thus, 
$i \partial(f_*\varphi) \wedge \bar \partial (f_*\varphi)$ has no mass on $X \backslash \Omega.$ It follows that 
$$i \partial(f_*\varphi) \wedge \bar \partial (f_*\varphi) \le d_k \bold{1}_\Omega[\ddc \big((f^{m_1+1})_*\psi\big)+ (f^{m_1+1})_*\eta] \le d_k[\ddc \big((f^{m_1+1})_*\psi\big)+ (f^{m_1+1})_*\eta]$$
because the last current is positive by Lemma \ref{le_quasipull}.  Combining this with (\ref{ineSP}) and (\ref{inechuanL1fsaoW12}) gives (\ref{inefsaofW12}).  The desired assertion then follows. The proof is finished. 
\endproof


Let $\varphi \in W^{1,2}_{*,f}$ and $\mathbf{m}=(m_0,m_1)$ a size representative of $\varphi.$   Consider  $f_*$ acting on Borel measurable functions.  Recall that $f_*$ preserves the set of constant functions.  As in the last section, let $b_0:= \int_X \varphi d \mu_0,$    and $\varphi_0:= \varphi -b_0.$  We define two sequences $\varphi_n, b_n$ as follows. Put
$$b_n= b_n(\varphi):= \int_X (f_*\varphi_{n-1}) d\mu_0,  \quad \varphi_n:= f_*\varphi_{n-1}- b_{n}$$
 for $n\ge 1.$ Note that $\varphi_n$ differs from $((f^n)_* \varphi)$ by a constant.   Lemma \ref{leppushvounded} yields that  $\mathbf{m}_n:=(d_k^n m_0,m_1+n)$ is a size representative of $\varphi_n.$ This coupled with Lemma \ref{le_repretrituyedoi} implies that 
 
\begin{lemma}\label{lesizerepesentaitive} $\mathbf{m}_n:=(d_k^n m_0,m_1+n)$ is also a size representative of $|\varphi_n|, \varphi_n^+$ and $\varphi^-_n.$ 
\end{lemma} 

 By Lemma \ref{le-sizereprese}, we get 
\begin{align} \label{chuan12varphidvarephi}
\| d \varphi_n \|_{L^2} \le  c_\epsilon m_0^{1/2}d_k^{n/2}(d_{k-1}+\epsilon)^{(n+ m_1)/2} 
\end{align}
 Using  (\ref{chuan12varphidvarephi}), (\ref{ineSP}) and  (\ref{inefsaofW12})   gives
 \begin{align} \label{chuan12varphidvarephi2}
 \| \varphi_n \|_{L^2} \le c_\epsilon m_0^{1/2} d_k^{n/2}(d_{k-1}+\epsilon)^{(n+ m_1)/2}, \quad |b_n| \le c_\epsilon m_0^{1/2} d_k^{n/2}(d_{k-1}+\epsilon)^{(n+ m_1)/2}
\end{align} 
for $n\ge 1$ and  some possible different constant $c_\epsilon.$   We are now in a situation very similar to that in the last section.  Using similar arguments as in the last section, we can show that $\lim_{n \to \infty} \langle d_k^{-n}(f^n)^* \omega^k, \varphi \rangle$ exists and denote by $b'_\infty(\varphi)$ this limit. Actually, we have
$$b'_\infty= \sum_{j=0}^\infty d_k^{-j} b_j.$$
It follows that 
\begin{align}\label{ine-chantrenbinfty}
|b'_\infty(\varphi)| \le   \| \varphi \|_{L^1}+  c_\epsilon  m_0^{1/2}  (d_{k-1}+\epsilon)^{m_1/2}
\end{align}
for some constant $c_\epsilon$ independent of $\varphi.$ Clearly, if $\varphi$ is a bounded quasi-p.s.h. function, $b'_\infty$ is equal to the same number defined in the last section.  Hence we have  $$\langle \mu_f, \varphi \rangle= b'_\infty(\varphi)$$ for bounded quasi-p.s.h. function $\varphi.$  Let $W^{1,2}_{**,f}$ be the subset of  $W^{1,2}_{*,f}$ consisting of functions which is continuous outside a closed pluripolar set. Observe that $f_*$ preserves $W^{1,2}_{**,f}$ because $f$ is a covering outside an analytic subset of $X$.   We now claim that 

\begin{lemma}\label{le-continuityoutsidepluripolar}   For  $\varphi \in W^{1,2}_{**,f},$ we have  $\langle \mu_f, \varphi \rangle= b'_\infty(\varphi).$ 
\end{lemma}

\proof  The proof is similar to that of \cite[Le. 5.5]{DS_decay}. We recall here for the readers' convenience.  We prove first that $\varphi$ is $\mu_f$-integrable.  We assume for the moment  that $\varphi\ge 0.$  Let $V$ be a closed pluripolar set such that $\varphi$ is continuous outside $V.$ Recall that $\mu_f$ has no mass on pluripolar sets, hence, on $V.$ Since $d_k^{-n}(f^n)^* \omega^k$ converges to $\mu_f$ as positive measures and $X \backslash V$ is open, we get 
$$\langle \mu_f, \varphi \rangle \le \liminf_{n \to \infty} \langle d_k^{-n}(f^n)^* \omega^k, \varphi \rangle=\lim_{n \to \infty} \sum_{j=0}^n d_k^{-j} b_j+ \liminf_{n \to \infty} \langle \omega^k, d_k^{-n}\varphi_n\rangle$$
which is equal to $b'_\infty(\varphi).$  Hence $\varphi$ is $\mu_f$-integrable if $\varphi \ge 0.$ In general, write $\varphi= \varphi^+ -\varphi^-$ and applying the last property shows that $\varphi$ is $\mu_f$-integrable.  If $\mathbf{m}=(m_0, m_1)$ is a size representative of $\varphi,$ then we also obtain that 
\begin{align} \label{estimatebprime}
|\langle \mu_f, \varphi | \le |b'_\infty(\varphi^+)|+ |b'_\infty(\varphi^-)| \le   c_\epsilon( \|\varphi\|_{L^1}+  m_0^{1/2} (d_{k-1}+\epsilon)^{m_1/2}),
\end{align}
for some constant $c$ independent of $\varphi.$   Now using $f^* \mu_f= d_k \mu_f$ gives
\begin{align*}
|\langle \mu_f, \varphi\rangle - b'_\infty(\varphi)|= |\langle \mu_f, d_k^{-n} (f^n)_* \varphi - b'_\infty(\varphi) \rangle| \le |c_n|+ |\langle \mu_f, d_k^{-n} \varphi_n \rangle|, 
\end{align*} 
where $c_n:= -\sum_{j\ge n+1}d_k^{-j} |b_j|.$ Observe that the first term in the right-hand side of the last inequality tends to $0$ because of (\ref{chuan12varphidvarephi2}). On the other hand, by (\ref{estimatebprime}) and Lemma \ref{lesizerepesentaitive},  the second term is bounded by 
$$c_\epsilon d_k^{-n}( \|\varphi_n\|_{L^1}+  m_0^{1/2} d_k^{n/2}(d_{k-1}+\epsilon)^{(m_1+n)/2})$$
which tends to $0$ as $n \to \infty.$ This yields the desired equality. The proof is finished.  
\endproof
 
 \begin{theorem} \label{thW12saosao} Let $X, f, d_k, d_{k-1}, \epsilon$ be as above with $d_k > d_{k-1}+\epsilon.$  Then there exists a constant $c_\epsilon$ such that 
 $$I_n(\psi, \varphi):= |\langle \mu_f, (\psi \circ f^n) \varphi\rangle - \langle \mu_f, \psi \rangle \langle \mu_f, \varphi \rangle| \le c_\epsilon \|\psi\|_\infty A_n(\varphi),$$
 where
$$A_n(\varphi):= \big[\|\varphi\|_{L^1}+m_0^{1/2}(d_{k-1}+\epsilon)^{m_1/2}\big] d_k^{-n/2} (d_{k-1}+\epsilon)^{n/2},$$
for every $\psi \in L^\infty(\mu_f)$, $\varphi \in W^{1,2}_{** ,f}$ and $(m_0, m_1)$ a size representative of $\varphi,$
 \end{theorem}

Note that if $\varphi$ is a bounded $\eta$-p.s.h. function for some continuous $(1,1)$-form $\eta$ of $L^\infty$-norm $\le 1,$ then there exists a constant $\tilde{m}_0$ independent of $\varphi$ such that $(\tilde{m}_0,1)$ is a size representative of $\varphi.$ Hence the above theorem gives a decay of correlation uniformly for every such $\varphi.$
 
\proof  Let the notations be as above.  Since $I_n(\psi, \varphi+c)= I_n(\psi, \varphi)$ for every constant $c$ because of the invariance of $\mu_f.$ We can assume that $\langle \mu_f, \varphi \rangle=0.$ By Lemma \ref{le-continuityoutsidepluripolar}, we get $b'_\infty(\varphi)=0.$ Hence,  $d_k^{-n}(f^n)_*(\varphi)= c_n+ d_k^{-n}\varphi_n.$    
Using $f^* \mu_f= d_k \mu_f$ gives 
\begin{align}\label{eqbiendoiIpsi}
I_n(\psi, \varphi)=  d_k^{-n} |\langle \mu_f, \psi  (f^n)_*(\varphi)\rangle|= |\langle \mu_f, \psi (c_n+ d_k^{-n}\varphi_n)\rangle|\le |c_n|+ d_k^{-n} |\langle \mu_f, |\varphi_n| \rangle|.
\end{align}
Note that as before we have 
$$|c_n| \le c_\epsilon A_n(\varphi)$$
for some constant $c_\epsilon$ independent of $\varphi.$
On the other hand, $f_*$ preserves $W^{1,2}_{** ,f},$ hence $\varphi_n \in W^{1,2}_{**,f}$ and so is $|\varphi_n|.$ By Lemma \ref{lesizerepesentaitive}, $(d_k^n m_0, m_1+n)$ is a size representative of $|\varphi_n|$ if $(m_0,m_1)$ is a size representative of $\varphi.$  Arguing as in the proof of Lemma \ref{le-continuityoutsidepluripolar} gives that 
$$d_k^{-n} |\langle \mu_f, |\varphi_n| \rangle | \le  c_\epsilon A_n(\varphi)$$
for some constant $c_\epsilon$ independent of $\varphi.$  Hence the desired inequality follows. This finishes the proof. 
\endproof
 
 \begin{proof}[End of the proof of Theorem \ref{the_mixing}]
 The central limit theorem for $\mu_f$ is a direct consequence of its decay of correlation as shown in \cite{DS_decay}. Hence it remains to prove the decay of correlation property. By Theorem \ref{thW12saosao}, for every $\mathcal{C}^1$ function $\varphi$ on $X,$ we have
 $$I(\psi, \varphi) \le c_\epsilon \|\psi\|_\infty \|\varphi\|_{\mathcal{C}^1} d_k^{-n/2} (d_{k-1}+\epsilon)^{n/2}.$$
This combined with the interpolation inequality for functionals on the Banach spaces $\mathcal{C}^1,\mathcal{C}^0$ gives the desired decay of correlation for $\mu_f,$ see \cite[p. 765]{DS_decay}. 

Recall that $\mu_f$ is K-mixing if for every $\varphi\in L^2(\mu_f),$ we have 
\begin{align}\label{ineKmixing}
\sup_{\psi\in L^2(\mu_f)} I_n(\psi,\varphi) \to 0.
\end{align}
Note that the operator $d_k^{-1} f_*$ can be extended to be a continuous linear operator  on $L^2(\mu_f)$ because $|f_* \varphi|^2 \le d_k f_*(|\varphi|^2)$.  As above, in order to prove (\ref{ineKmixing}), we can assume that $\langle \mu_f, \varphi\rangle=0.$ Using (\ref{eqbiendoiIpsi}) gives 
\begin{align}\label{inevhanIpsivarphiL2}
I(\psi, \varphi) \le  \|d_k^{-n} (f^n)_* \varphi\|_{L^2(\mu_f)}.
\end{align}
Consider now $\varphi$ to be a bounded function in $W^{1,2}_{**,f}$. The set of these functions is dense in $L^{2}(\mu_f).$ We have 
$$ \|d_k^{-n} (f^n)_* \varphi\|_{L^2(\mu_f)} \le  \|\varphi\|_\infty  \|d_k^{-n} (f^n)_* \varphi\|_{L^1(\mu_f)}$$
which tends to $0$ by the proof of Theorem \ref{thW12saosao}. This combined with (\ref{inevhanIpsivarphiL2}) gives (\ref{ineKmixing}).  The proof is finished. 
 \end{proof}

\begin{remark} By the inequality (\ref{chuan12varphidvarephi2}),  we see that for every complex measure $\nu$ with $L^2$ density and $\nu(X)=1,$ $d_k^{-n} (f^n)^* \nu$ converges weakly to $\mu_f.$
\end{remark}

\section{Examples}  \label{example}

In this section, we present examples of holomorphic dynamical systems on non-K\"ahler manifolds. Some of them were already considered by Gromov \cite{Gromov_entropy}. 

\subsection{Hopf manifolds.}  Let $H$ be the standard Hopf manifold $(\C^k\backslash \{0\})/ \{z \sim \lambda  z\},$ for $\lambda \in \C^*$ and $k \ge 2.$   Recall that the natural map $p_H: H \to \P^{k-1}$ defined by $(z_1,\ldots, z_k) \mapsto [z_1: \cdots :z_k]$ is a fiber bundle whose fibers are $\C^*/\{t \sim \lambda t\}$ which are  compact Riemann surfaces of genus $2.$ Every  holomorphic endomorphism $f$ of $H$ is induced by an endomorphism $F$ of $\C^k\backslash \{0\}$ whose components are homogeneous polynomials of the same degree $\tilde{d}.$ Hence $f$ is open and of finite fibers.  Clearly $F$ induces naturally a holomorphic endomorphism $f'$ of $\P^{k-1}$ and 
\begin{align}\label{eq_giaphoandff}
p_H \circ f= f' \circ p_H.
\end{align}
Using this we get  $$d_k(f)= \tilde{d}^{k+1}, \quad P_n= \tilde{d}^{k+1}$$
where $P_n$ is the set of periodic points of period $n$ of $f.$ By \cite{Gromov_entropy}, the topological entropy $h_t(f)$ is equal to  $\log d_k(f).$

Let $r$ be a strictly positive number. Let $\cali{D}(r, \tilde{d})$ be the set of holomorphic maps $F: \C^k \backslash \{0\} \to \C^k \backslash \{0\}$ whose components are homogeneous polynomials of degree $\tilde{d}$ for which 
\begin{align} \label{ine_dinhghiaDngar}
\|DF(z)\| \le r \tilde{d} \|F(z)\|/ \|z\|,
\end{align} 
for $z \in \C^k \backslash \{0\},$ where $D F$ is the differential of $F$ and   the norms are the Euclidean norms. Using $$\|D F^n(x)\| \le \|DF\big(F^{n-1}(x)\big)\| \cdots \|DF\big(F(x)\big)\| \cdot \| DF(x) \|$$
 and  (\ref{ine_dinhghiaDngar}), we infer that 
\begin{align} \label{ine_dinhghiaDngar2}
\|DF^n(z)\| \le r^n \tilde{d}^n \|F^n(z)\|/ \|z\|,
\end{align} 
for $z \in \C^k \backslash \{0\}$ and $n \in \N^*.$   We can see that the map $F_0(z):= (z_1^{\tilde{d}}, \ldots, z_k^{\tilde{d}})$ belongs to $\cali{D}(2k, \tilde{d})$ because 
$$\|D F_0(z)\|= \tilde{d}\big(\sum_{j=1}^k |z_j|^{2(\tilde{d} -1)}\big)^{1/2} \le  2k \tilde{d} \big(\sum_{j=1}^k |z_j|^{2\tilde{d}}\big)^{1/2}\big(\sum_{j=1}^k |z_j|^{2}\big)^{-1/2}= 2 k \tilde{d} \|F_0(z)\| / \|z\|.$$
We can construct easily some other examples. 
     
\begin{lemma} \label{le_ex_dominanthopf} Let $\tilde{d}$ be  a positive integer $\ge 2$ and $r$ a positive real number $\ge 1$.   Let $f: H  \to H$ be the holomorphic map induced by a map $F \in \cali{D}(r, \tilde{d}).$ 
 Then we have  $ d_q(f) \le  r^2 \tilde{d}^{q+1} $ for every $0 \le q \le k-1.$ In particular, if $\tilde{d}> r^2,$ then  $f$ has a dominant topological degree, \emph{i.e,} $d_k>d_q$ for $0 \le q \le k-1$. 
\end{lemma}

Arguing as in \cite[Pro. 2.7]{DS_book},  we can see that as in the case of polynomial-like maps, the property of having a dominant topological degree is preserved under  small perturbations. Hence the above lemma provides us a rich class of self-maps  with a dominant topological degree. 

\proof Let the notation be as above. We already know that $d_k= \tilde{d}^{k+1}.$ Hence, if we can prove the desired inequality, then  $f$ has a dominant topological degree provided that $\tilde{d}>(r^2+2).$ Recall that $d_0(f)=1.$  

Let $1 \le  q \le k-1.$ Let $\omega_\FS$ be the Fubini-Study form on $\P^{k-1}$ and $p_{\C^{k}}$ the natural projection from $\C^k \backslash \{0\}$ to $\P^{k-1}.$ The pull-back $p_{\C^{k}}^* \omega_{\FS}$ on $\C^{k} \backslash \{0\}$ is given by 
$$2 \,p_{\C^{k}}^* \omega_\FS=  \ddc\log(|z_1|^2+ \cdots+|z_k|^2).$$
Recall here $\ddc= i \partial \overline \partial/\pi.$ Define 
$$\omega:=(2\pi)^{-1} \|z\|^{-2} \sum_{j=1}^k dz_j \wedge d\bar z_j$$
and 
$$\eta:= \sum_{j=1}^k z_j d z_j, \quad \omega':= (2\pi)^{-1} \|z\|^{-4}\eta \wedge \overline \eta \ge 0.$$
Both $\omega, \omega'$ induce well-defined smooth forms on $H$ which are denoted by the same notations $\omega, \omega'$ respectively for simplicity.  Observe that $\omega$ is a Hermitian metric on $H.$   Direct computations show that 
\begin{align} \label{eq_omegFSphay}
\omega= p_{\C^{k}}^* \omega_\FS+ \omega', \quad  \omega' \wedge \omega' =0.
\end{align} 
By (\ref{eq_omegFSphay}), on $H,$ we have 
$$\omega^q= p_H^*\omega_{\FS}^q+ p_H^*\omega_{\FS}^{q-1} \wedge \omega'$$
for $1 \le q \le k-1.$  Using the last equality and (\ref{eq_giaphoandff}) gives
$$f^* \omega^q= p_H^* f'^* \omega_{\FS}^q+ p_H^*f'^*\omega_{\FS}^{q-1} \wedge f^*\omega'.$$
 It follows that
\begin{align} \label{eq_tinhd_qfHopf}
\int_H (f^n)^* \omega^q \wedge \omega^{k-q}& = \int_H p_H^*\big[ (f'^n)^* \omega_{\FS}^q\wedge \omega_{\FS}^{k-q-1}\big] \wedge  \omega'+   \int_H p_H^*\big[ (f'^n)^* \omega_{\FS}^{q-1}\wedge \omega_{\FS}^{k-q}\big] \wedge (f^n)^* \omega'\\
\nonumber
&+\int_H p_H^*\big[ (f'^n)^* \omega_{\FS}^{q-1}\wedge \omega_{\FS}^{k-q-1}\big] \wedge (f^n)^* \omega' \wedge \omega'.
\end{align}
Denote by $I_1, I_2,I_3$ the first, second  and third integrals in the right-hand side of the last equality. Using Fubini's theorem gives
\begin{align*}
I_1 =  \int_{[z] \in \P^{k-1}}(f'^n)^* \omega_{\FS}^q\wedge \omega_{\FS}^{k-q-1} \int_{p_H^{-1}([z])}  \omega'
\end{align*}
and
\begin{align*}
I_2 =\int_{[z] \in \P^{k-1}}(f'^n)^* \omega_{\FS}^{q-1} \wedge \omega_{\FS}^{k-q} \int_{p_H^{-1}([z])} (f'^n)^*\omega'=\tilde{d}^{2n} \int_{[z] \in \P^{k-1}}(f'^n)^* \omega_{\FS}^{q-1} \wedge \omega_{\FS}^{k-q} \int_{p_H^{-1}([z])} \omega'
\end{align*}
because the topological degree of $f_{[z]}$ is equal to $\tilde{d}^2.$ Thus we get
\begin{align}\label{eq_tinhI12}
\lim_{n\to \infty} I_1^{1/n}= d_{q}(f')= \tilde{d}^{q}, \quad  \lim_{n\to \infty} I_2^{1/n}= \tilde{d}^2 \,  d_{q-1}(f')= \tilde{d}^{q+1},
\end{align} 
where recall that the dynamical degree $d_q(f')$ of $f'$ is equal to $\tilde{d}^{q}$ for $1 \le q\le k-1.$ 

It remains to estimate $I_3.$ 
By (\ref{eq_omegFSphay}), 
$$(f^n)^* \omega' \wedge \omega'= (f^n)^* \omega \wedge \omega'- p_H^* (f'^n)^* \omega_\FS \wedge \omega'.$$
Direct computations show that
$$(F^n)^* \omega= (2\pi)^{-1}\|F^n\|^{-2} \sum_{j=1}^k d F^n_j \wedge d \overline{F^n}_j \le  (2\pi)^{-1} k (r\tilde{d})^{2n} \|z\|^{-2} \sum_{j=1}^k  d z_j \wedge d \overline{z}_j=k (r\tilde{d})^{2n} \omega $$
by  (\ref{ine_dinhghiaDngar2}). It follows that 
$$ (f^n)^* \omega' \wedge \omega' \le  k (r\tilde{d})^{2n}  \omega \wedge \omega'+ p_H^* (f'^n)^* \omega_\FS \wedge \omega' \le k [(r\tilde{d})^{2n}+1]  p_H^* \omega_\FS \wedge \omega'$$
which implies 
$$I_3 \le  k [(r\tilde{d})^{2n}+1] \int_H p_H^*\big[ (f'^n)^* \omega_{\FS}^{q-1}\wedge \omega_{\FS}^{k-q}\big] \wedge \omega'.$$
Taking the power $1/n$ in the last inequality and letting $n \to \infty$ give
\begin{align}\label{eq_tinhI123}
\limsup_{n\to \infty} I_3^{1/n} \le r^2  \tilde{d}^2 d_{q-1}(f')= r^2 \tilde{d}^{q+1}
\end{align}
for $1 \le q\le k-1.$ Combining (\ref{eq_tinhI123}), (\ref{eq_tinhI12}) and (\ref{eq_tinhd_qfHopf}) yields
$$d_q(f) \le \limsup_{n\to \infty}(I_1+ I_2+ I_3)^{1/n} \le  r^2  \tilde{d}^{q+1}.$$ 
The proof is finished. 
\endproof

\subsection{Calabi-Eckmann manifolds.}  Let $\alpha \in \C \backslash \R$. Let $k,l$ be integers $\ge 2.$  The Calabi-Eckmann manifold $X$ is defined by 
$$X:= (\C^k\backslash \{0\}) \times (\C^l\backslash \{0\})/ \sim,$$
where the equivalence relation $\sim$ is given by  $(z,w) \sim (e^t z, e^{\alpha t} w)$ for every $t \in \C.$ Recall that $X$ is diffeomorphic to $\S^{2k-1} \times \S^{2l-1}$ and non-K\"ahler because $H^2(X)=\{0\}.$  For every homogeneous polynomials $F(z),G(w)$ of the same degree,  the self-map of $ (\C^k\backslash \{0\}) \times (\C^l \backslash \{0\})$ given by $(z,w) \mapsto \big(F(z),G(w)\big)$ can descend to a self-map of $X.$ It is likely that we can obtain a good class of self-maps of $X$ with dominant topological degree as above.

\subsection{Nilmanifolds.}   Consider $G$ a complex Lie group and $\Gamma$ a closed complex Lie subgroup of $G$ such that $X:= G/\Gamma$ is a compact non-K\"ahler manifold of dimension $k$.  By \cite{Benson_Gordon}, nilmanifolds which are not tori  are examples of such manifolds. For such $X,$ every $g \in G$ and every  $A \in \Aut(G)$ preserving $\Gamma,$ the affine transformation $g A$ induces a holomorphic automorphism on $X.$ In the real setting, the dynamical systems associated to such maps possess interesting properies and has been studied extensively. We refer to \cite{Parry_nil,Gorodnik_Spatzier,book_DSEA} for informations.

\subsection{Blowups.} Let $X$ be a compact complex manifold. Let $\widehat X$ be a compact manifold bimeromorphic to $X$ via a map $\sigma: \widehat X \to X.$  Then given a meromorphic self-correspondence $f$ on $X$, $f_\sigma:=\sigma^{-1} \circ f \circ \sigma$ is a self-correspondence on $\widehat X.$ We can take, for example,  $X$ to be a nilmanifold, a Hopf manifold or a Calabi-Eckmann manifold and $\widehat X$ to be the blowup of $X$ along a smooth submanifold  $V$ of $X$ (a point for example).  By a well-known example of Hironaka \cite{Hironaka_ex}, there exist  a compact K\"ahler manifold $X$ and a non-K\"ahler manifold $\widehat X$ bimeromorphic to $X.$ Such manifolds $\widehat X$ are  in the class $\mathcal{C}$ of Fujiki.  Due to the lack of a K\"ahler form, we don't know whether the dynamical degrees of $f$ and $f_\sigma$ are the same even for $X$ K\"ahler. 

\bibliography{biblio_family_MA,bib_expose_cnrs,biblio_Viet_papers}
\bibliographystyle{siam}

\end{document}